\documentclass[letterpaper,10pt]{amsart}

\usepackage[all]{xy}                        
\CompileMatrices                            
\UseTips                                    

\usepackage[bookmarks=true]{hyperref}       
\usepackage{amssymb}                        
\usepackage{xspace}                         
\usepackage{epsfig}
\usepackage{stmaryrd}                        

\usepackage{layout}                        
\reversemarginpar

\theoremstyle{plain}
\newtheorem{theorem}{Theorem}
\newtheorem*{theorem*}{Theorem}

\newtheorem*{corollary}{Corollary}
\newtheorem{lemma}{Lemma}

\theoremstyle{definition}
\newtheorem{definition}{Definition}

\theoremstyle{remark}

\newcommand{\enm}[1]{\ensuremath{#1}}          
\newcommand{\op}[1]{\operatorname{#1}}
\newcommand{\cal}[1]{\mathcal{#1}}

\renewcommand{\bar}[1]{\overline{#1}}

\newcommand{\NN}{\enm{\mathbb{N}}}
\newcommand{\RR}{\enm{\mathbb{R}}}
\newcommand{\QQ}{\enm{\mathbb{Q}}}
\newcommand{\ZZ}{\enm{\mathbb{Z}}}

\renewcommand{\AA}{\enm{\mathbb{A}}}

\newcommand{\Jj}{\enm{\cal{J}}}

\newcommand{\Oo}{\enm{\cal{O}}}

\renewcommand{\phi}{\varphi}        
\renewcommand{\theta}{\vartheta}
\renewcommand{\epsilon}{\varepsilon}

\newcommand{\Ann}{\op{Ann}}         

\newcommand{\Spec}{\op{Spec}}

\newcommand{\Supp}{\op{Supp}}

\newcommand{\tensor}{\otimes}         

\renewcommand{\to}[1][]{\xrightarrow{\ #1\ }}

\newcommand{\defeq}{\stackrel{\scriptscriptstyle \op{def}}{=}}

\renewcommand{\aa}{\mathfrak{a}}

\newcommand{\mm}{\mathfrak{m}}
\newcommand{\ord}{\op{ord}}
\renewcommand{\div}{\op{div}}
\newcommand{\floor}[1]{\llcorner #1 \lrcorner}
\newcommand{\roof}[1]{\ulcorner #1 \urcorner}
\newcommand{\ceil}[1]{\ulcorner #1 \urcorner}

\newcommand{\Int}{\op{Int}}
\newcommand{\smallfrac}[2]{\textstyle{\frac{#1}{#2}}}

\newcommand{\Newt}{\op{Newt}}

\begin{document}

\title{Multiplier ideals and modules on toric varieties}
\author{Manuel Blickle}
\address{Universit\"at Essen, FB6 Mathematik, 45117 Essen, Germany}
\email{manuel.blickle@uni-essen.de}\urladdr{\url{www.mabli.org}}
\date{October 29, 2003}
\thanks{The author is grateful to Nobuo Hara for interesting discussions and thanks the referee
for a careful reading and thoughtful comments.}

\begin{abstract}
A formula computing the multiplier ideal of a monomial ideal on an arbitrary
affine toric variety is given. Variants for the multiplier module and test
ideals are also treated.
\end{abstract}

\maketitle


\section{Introduction and results}
In this note I generalize Howald's formula \cite{Howald.MultMon} for the
multiplier ideal of monomial ideals in a polynomial ring to the case of torus
invariant ideals in an arbitrary (normal) toric variety. If $\aa \subseteq
k[x_1,\ldots,x_n]$ is a monomial ideal, his formula computes the multiplier
ideal:
\[
    \Jj(\aa^c) = \langle x^v | v+(1,\ldots,1) \in \text{ interior of } c\Newt(\aa) \rangle
\]
Here and henceforth, $x^v$ is shorthand for $x_1^{v_1}\cdot\ldots\cdot
x_n^{v_n}$ which is a monomial in $k[x_1,\ldots,x_n]$. To each monomial $x^v$
one associates its exponent vector $v=(v_1,\ldots,v_n)$ inside the lattice of
exponents $M \cong \ZZ^n$. In this way we assign to a monomial ideal $\aa$ the
set of exponents of monomials in $\aa$; the convex hull of this set in $M_\RR=M
\tensor_\ZZ \RR$ is the \emph{Newton Polyhedron}, $\Newt(\aa)$, of $\aa$.

\subsection{Multiplier ideal on a pair}
The generalization of this formula involves the multiplier ideal of a singular,
not necessarily $\QQ$-Gorenstein, variety. Let us recall the setup, for details
consult \cite{PAG}, Chapter 9. Let $(X,\Delta)$ be a pair, consisting of a
normal variety $X$ and a $\QQ$-divisor $\Delta$ such that $K_X+\Delta$ is
$\QQ$-Cartier. This is the situation in which multiplier ideals can be defined.
\begin{definition}\label{def.multid}
    For such pair $(X,\Delta)$ and an ideal sheaf $\aa$ on $X$,
    choose a log resolution $\mu: Y \to X$  of $\aa$ which at the same time is a log
    resolution of the pair. In particular, $\aa\cdot\Oo_Y=\Oo_Y(-A)$ for some
    effective normal crossing divisor $A$. Then for all $c>0$, the \emph{multiplier ideal}
    of $\aa^c$ on the pair $(X,\Delta)$ is
    \[
        \Jj((X,\Delta),\aa^c) = \mu_*\Oo_Y(K_Y - \floor{\mu^*(K_X+\Delta) +
        cA}).
    \]
\end{definition}
The expression on the right does \emph{not} make sense a priori, since it
depends on the divisor chosen to represent the canonical class (rounding does
not commute with numerical equivalence). But this definition can be made
rigorous using discrepancies as explained in \cite[Chapter 9.3.F]{PAG}. In the
case of toric varieties there is a canonical torus invariant representative
$K_X$ of the canonical class. Namely, if $D_1,\ldots,D_s$ are the torus
invariant prime Weil divisors on $X$ then $K_X=-(D_1+\ldots+D_s)$. Using this
choice for $K_X$ and $K_Y$, one easily verifies that $K_Y-\mu^*(K_X+\Delta)$ is
equal (not just numerically equivalent) to $E-\widetilde{\Delta}$, where $E$ is
$\mu$-exceptional and $\widetilde{\Delta}$ denotes the strict transform of
$\Delta$. This validates our setup.

Now let $(X,\Delta)$ be a pair such that $X$ is an (affine) toric variety (say
$X = \Spec R$ for some normal semigroup ring $A \subseteq k[x_1^{\pm
1},\ldots,x_n^{\pm 1}]$) and $\Delta$ is a torus invariant $\QQ$-divisor. Since
$K_X+\Delta$ is $\QQ$-Cartier and torus invariant, there is a monomial $x^u$
such that $\div x^u = r(K_X+\Delta)$ for some integer $r$. With $w=u/r$, one
obtains the following result:
\begin{theorem}\label{thm.ideal}
    Let $\aa$ be a monomial ideal on $X$. Then
    \[
        \Jj((X,\Delta),\aa^c) = \langle x^v \in R| v + w \in \text{ interior of } c\Newt(\aa)
        \rangle
    \]
    for all $c>0$.
\end{theorem}
If $X=\AA^n$ and $\Delta=\emptyset$ this is just Howald's forlmula. If $X$ is
$\QQ$-Gorenstein and $\Delta = \emptyset$ this recovers a generalization by
Hara and Yoshida \cite{HaraYosh} obtained with positive characteristic methods;
also Howard Thompson obtained this case independently.

\subsection{Multiplier module}
A variant of the multiplier ideal on a singular variety is the
multiplier module. It has the advantage to be defined on any normal variety,
not relying of a $\QQ$-Gorenstein assumption, see, for example,
\cite{HySmi.Kawamata}. The multiplier module is no longer an ideal but a
submodule of the canonical module $\omega_X$.
\begin{definition}
    Let $X$ be a normal variety and let $\aa$ be a sheaf of ideals on $X$. Let $\mu: Y
    \to X$ be a log resolution of $\aa$ such that $\aa \cdot \Oo_Y =
    \Oo_Y(-A)$. Then the \emph{multiplier module} is defined as
    \[
        \Jj_{\omega}(\aa^c) = \mu_*\Oo_Y(K_Y-\floor{cA}) \subseteq \omega_X
    \]
    for all $c>0$.
\end{definition}
This is independent of the chosen log resolution of $\aa$. For an affine toric
variety $X$ there is a canonical embedding of $\omega_X \subseteq \Oo_X$, which
realizes $\omega_X$ as a monomial ideal. In this canonical description of
$\omega_X$ the formula computing the multiplier module takes a stunningly easy
form.
\begin{theorem}\label{thm.module}
    Let $X$ be an affine toric variety and $\aa$ a monomial ideal. Then
    \[
        \Jj_\omega(\aa^c) = \langle x^v | v \in \text{ interior of } c
        \Newt(\aa)\rangle \subseteq \omega_X
    \]
    for all $c>0$.
\end{theorem}
Clearly, as multiplier ideals are local in nature one can derive similar
statements for the multiplier ideal or module of an arbitrary normal toric
variety by computing on an affine toric cover using the above formulas. In
general the notion of multiplier module and ideal require the characteristic to
be zero since their construction employs resolution of singularities. Since for
toric varieties log resolutions are available also in positive characteristic,
the preceding results are valid, independent of the characteristic.

\subsection{Positive characteristic: Test ideals}
In positive characteristic, Hara and Yoshida introduced a related notion of the
test ideal $\tau(\aa^c)$ which does not require a $\QQ$-Gorenstein assumption.
They also show that in the $\QQ$-Gorenstein case the multiplier ideal
$\Jj(\aa^c)$ generically yields their test ideal $\tau(\aa^c)$ under reduction
to positive characteristic. As a final result I give a formula (see Theorem
\ref{thm.poschar}) for the test ideal $\tau(\aa^c)$ for an affine toric variety
defined over a field of positive characteristic. This formula specializes to
the formula above in the $\QQ$-Gorenstein case. We relegate the statement
(Theorem \ref{thm.poschar}) and proof to the last section since we first need
to recall some basics from toric geometry.

\section{Toric Setup}
Following Fulton \cite{Fulton.Toric}, I fix a dual pair of lattices $N = M^\vee
\cong \ZZ^n$. Let $\sigma \subseteq N_\RR=N\tensor_\ZZ \RR$ be a strongly
convex rational polyhedral cone given by $\sigma=\{\,r_1u_1+\ldots+r_tu_t\, |\,
r_i \in \RR_+ \,\}$  for some $u_1,\ldots,u_t \in N$. The dual cone
$\sigma^\vee$ is a (rational convex polyhedral) cone in $M_\RR$ defined by
$\sigma^\vee = \{\, m \in M_\RR\, |\, (m,v) \geq 0,\, \forall\, v\in \sigma
\,\}$ where $(\cdot,\cdot)$ denotes the pairing of the dual lattices $M$ and
$N$. The lattice points in $\sigma^\vee$ give rise to a sub-semigroup of
Laurent polynomials $k[x_1^{\pm 1},\ldots,x_n^{\pm 1}]$ generated by those
monomials $x^m=x_1^{m_1}\cdot\ldots\cdot x_n^{m_n}$ such that $m \in
\sigma^\vee$. Identifying a monomial with its exponent this yields the affine
semigroup ring
\[
    R_\sigma = k[ \sigma^\vee \cap M ].
\]
The affine toric variety is $X_\sigma = \Spec R_\sigma$. Since $R_\sigma$ is
contained in the ring of Laurent polynomials, $X_\sigma$ contains the torus
$(k^*)^n = T^n=\Spec k[x_1^{\pm 1},\ldots,x_n^{\pm 1}]$ as a dense subset. The
action of the torus on itself extends naturally to an action on the whole of
$X_\sigma$.

An ideal $\aa$ of a toric ring $R_\sigma$ is invariant under the action of the
torus if and only if it is generated by monomials, that is, if it is a monomial
ideal. To such ideal $\aa$ one associates its \emph{Newton polyhedron}
$\Newt(\aa) \subseteq M_\RR$, defined as the convex hull of $\aa$, or more
precisely the convex hull of the set of exponents $m$ of the monomials $x^m$ in
$\aa$.

\subsection{Divisors}\label{sec.divisors} The prime (Weil) divisors which are fixed by this torus action have
an easy description, they correspond to the edges (or rays) of the cone
$\sigma$. Let $v_1,\ldots,v_s$ be the first lattice points on the edges of
$\sigma$. Their orthogonals $v_i^\perp \cap \sigma^\vee \subseteq \sigma^\vee$
are the facets of $\sigma^\vee$. They define codimension one subvarieties $D_i$
of $X_\sigma$. Any torus invariant divisor can be written as a sum of the $D_i$
and we denote the lattice of torus invariant divisors by $L=L^{X_\sigma}=\oplus
\ZZ D_i$. A torus invariant Cartier divisor can be written as $\div x^m$ for a
Laurent monomial $x^m$. The fact that $\ord_{D_i} x^m = (m,v_i)$ shows
therefore that a divisor $D=\sum d_iD_i$ is $\QQ$-Cartier if and only if there
is $w \in M_\QQ$ such that $(w,v_i)=d_i$ for all $i$.

There is a canonical choice of a divisor $K_{X_{\sigma}}$ to represent the
canonical class. Namely the torus invariant divisor which is just the negative
of the sum of the prime divisors $D_i$, that is $K_{X_\sigma} = -\sum D_i$.
With this canonical divisor, the canonical module $\omega_{X_\sigma}$ is the
ideal of $R_\sigma$ consisting precisely of the monomials $x^m$ such that $m$
is in the interior of $\sigma^\vee$. In other words, $x^m \in
\omega_{X_\sigma}$ if and only if $\ord_{D_i} x^m = (m,v_i) > 0$ for all $i$.

A general toric variety is made up from affine pieces via the datum of a fan
$\Sigma$. A fan is a collection of compatible cones (any two cones in $\Sigma$
meet in a common face \ldots). The considerations made above are also valid for
this non affine setting, refer to \cite{Fulton.Toric} or \cite{Danilov.Toric}
for justification of all of the above (and all that follows in this section).

\subsection{Resolutions} A toric variety $X$ can be desingularized via an easy
combinatorial procedure involving subdividing the fan $\Sigma$ to arrive at a
fan $\Sigma'$ all of whose maximal cones are spanned by a basis of $N$. We only
need to know that this yields a torus equivariant desingularization $\mu:Y \to
X$ where $Y$ is the toric variety defined by $\Sigma'$. Similarly, log
resolutions of an ideal are also obtained torically. Since the monomial ideals
are precisely the torus invariant ideals, it follows that the multiplier ideals
of monomial ideals will also be monomial ideals.

Recall that toric varieties have at worst rational singularities, therefore,
$\mu_* \omega_{Y} = \omega_X$. This has the following immediate consequence
which will be used later.
\begin{lemma}\label{lem.pos}
    Let $\mu: Y \to X$ be a toric resolution of singularities. If $x^m \in \omega_X$
    then $\ord_{D_i} \mu^* x^m > 0$ for $D_i$ any torus invariant prime Weil divisor
    on $Y$. Consequently, $\Supp (\mu^* \div x^m)$ and $K_Y$ have the same support.
\end{lemma}
\begin{proof}
    The condition that $x^m$ is a section of $\omega_X$ implies by $\mu_* \omega_{Y} = \omega_X$
    that $x^m$ (viewed on $Y$) is a section of $\omega_Y$. But this just means
    $\ord_{D_i}\mu^*x^m > 0$ for all $i$ as claimed.
\end{proof}
Another simple observation is the following lemma.
\begin{lemma}\label{lem.formula}
    Let $\mu: Y \to X$ be a toric log resolution of the monomial ideal $\aa$
    such that $\aa\cdot\Oo_Y=\Oo_Y(-A)$. For a monomial $x^m$ one has
    \[
        cm \in c'\Newt(\aa) \iff c\mu^* \div x^m \geq c'A.
    \]
    for all rational (real) numbers $c$ and $c'$.
\end{lemma}
\begin{proof}
    It is well known \cite{Fulton.Toric} that $m \in
    \Newt(\aa)$ if and only if $x^m \in \bar{\aa}$, the integral closure of
    $\aa$. Since $\bar{\aa}=\mu_*\Oo_Y(-A)$, this is equivalent to $\mu^*\div x^v \geq A$. This
    was the case $c=c'=1$ and the rest follows easily by expressing $c$ and $c'$ as
    integer fractions and clearing denominators.
\end{proof}
With these preparations we proceed to the proofs of Theorem \ref{thm.ideal} and
Theorem \ref{thm.module}.
\section{Proofs of Theorem 1 \& 2}
Both proofs are essentially the same argument. I first give the proof for
Theorem \ref{thm.module} and then indicate the small changes needed in the
proof of Theorem \ref{thm.ideal}.
\begin{proof}[Proof of Theorem \ref{thm.module}]
    Fix a toric log resolution $\mu: Y \to X_\sigma$ of $\aa$ such that $\aa\cdot
    \Oo_Y = \Oo_Y(-A)$. For $m$ to be in the interior of the Newton polyhedron
    $c\Newt(\aa)$ is the same as
    \[
        m - \epsilon m' \in c\Newt(\aa)
    \]
    for all (some) $m' \in M$ in the interior of $\sigma^\vee$ and all small enough
    (rational) $\epsilon$. This is because if $z \in \Newt(\aa)$, then
    $z+\sigma^\vee \subseteq \Newt(a)$. Applying Lemma \ref{lem.formula} we get that this is
    equivalent to
    \[
        \mu^* \div x^m - \epsilon \mu^* \div x^{m'} \geq cA,
    \]
    which is, since $\mu^* \div x^m$ is integral, equivalent to
    \[
        \mu^* \div x^m \geq \roof{\epsilon \mu^* \div x^{m'}+cA}.
    \]
    By Lemma \ref{lem.pos} (using that $m'$ is in
    the interior of $\sigma^\vee$ if and only if $x^{m'}\in \omega_X$),
    the $\QQ$-divisor $\epsilon \mu^* \div x^{m'}$ is effective with support
    equal to $-K_Y$, that is the prime divisors in $\mu^* \div x^{m'}$ are precisely
    all the torus invariant divisors. Then an exercise in rounding shows that the last
    inequality is equivalent to
    \[
        \mu^*\div x^v \geq K_Y + \floor{cA}
    \]
    as claimed.
\end{proof}

\begin{proof}[Proof of Theorem \ref{thm.ideal}]
    Let $\mu:Y \to X$ be a toric log resolution as in the last proof. The condition that
    $m+w$ is in the interior of the Newton polyhedron $\Newt(\aa)$
    is equivalent to $m+u/r-\epsilon m' \in \Newt(\aa)$ for $m'$ in the interior
    of $\sigma$ and small enough rational $\epsilon > 0$. Again, by Lemma
    \ref{lem.formula} this is equivalent to
    \[
        \mu^*\div x^m + \frac{1}{r}\mu^*\div x^u - \epsilon \mu^*\div x^{m'} \geq cA.
    \]
    Using that $\div x^u = r(K_X+\Delta)$ and proceeding analogous to the last
    proof this is equivalent to    \[
        \mu^* \div x^m \geq K_Y -\floor{\mu^*(K_X+\Delta)+cA}
    \]
    which says, by definition, nothing but that $x^m \in \Jj((X,\Delta);\aa^c)$.
\end{proof}

\section{Formula for the test ideal and speculations}
In positive characteristic, Hara and Yoshida \cite{HaraYosh} introduced a
notion of \emph{test ideal} which, under the process of reduction to positive
characteristic corresponds to the multiplier ideal. This notion grew out of
tight closure theory \cite{Huneke.TightBook}, where test ideals play an
important role. Their definition does not rely on a $\QQ$-Gorenstein
assumption, which can be attributed to the fact that their construction does
not use resolutions of singularities. With their test ideal $\tau(\aa^c)$
replacing the multiplier ideal one obtains a result for the test ideal of a
monomial ideal on an arbitrary affine normal toric variety over a field of
positive characteristic.
\begin{theorem}\label{thm.poschar}
    Let $X_\sigma=\Spec R$ be a toric variety over a field of positive characteristic and $\aa$
    a monomial ideal. Then a monomial $x^m \in R$ is in $\tau(\aa^c)$ if and only if
    there exists $w \in M_\RR$ with $(w,v_i)\leq 1$ for all $i$, such that
    \[
        m+w \in \text{ interior of }c\Newt(\aa).
    \]
    If $R$ is $\QQ$-Gorenstein, then there is a $w_0$ with $(w_0,v_i)=1$ for all $i$.
    Therefore $x^m \in \tau(\aa^c)$ if and only if $m+w_0 \in \text{ interior of }c\Newt(\aa)$.
\end{theorem}
The definition of the test ideal associated to an ideal $\aa \subseteq R$ and a
rational parameter $c$ is
\[
    \tau(\aa^c) \defeq \{\, h \in R\, |\, hI^{*\aa^c} \subseteq I\text{ for all ideals $I$ of $R$ }\}.
\]
The $\aa^c$-tight closure appearing is defined similarly as the usual tight
closure: $x \in I^{*a^c}$ if there is an $h \in R^\circ$ such that for all
$q=p^e$ one has $hx^q \aa^{\ceil{cq}} \subseteq I^{[q]}$. I will not elaborate
on the properties of $\aa^c$--tight closure and the resulting test ideals (see
\cite{HaraYosh}) but instead mention the one result needed to prove Theorem
\ref{thm.poschar}.

\begin{lemma}\label{lem.testIdeal}
    Let $(R,\mm)$ be a graded ring and $\aa$ a homogeneous ideal. Then
    \[
        \tau(\aa^c) = \Ann_R (0^{*\aa^c}_{E_{R/\mm}}),
    \]
    where $E_{R/\mm}$ is the injective hull of the residue field of $R$.
\end{lemma}
\begin{proof}
First observe that Theorem 3.3 of \cite{LyuSmith.StrongWeak} (saying that in
the graded case and for an Artinian module the \emph{finitistic tight closure}
(which is a variant of tight closure) is equal to the tight closure) can easily
be adapted to the case of $\aa^c$--tight closure, provided that $\aa$ is also
graded. Since semigroup rings are in particular $\NN$--graded and so are all
our modules and ideals involved, this implies that
\[
        0^{*\aa^c}_E = \bigcup_{\substack{\text{finitely}\\ \text{generated}\\ V \subseteq
        E}}
        0^{*\aa^c}_V.
\]
As by \cite{HaraYosh}, Proposition 1.9, the test ideal $\tau(\aa^c)$ is the
annihilator of the right hand side the result follows.
\end{proof}

With this observation we can modify the proof of Theorem 4.8 of \cite{HaraYosh}
to apply also to the non $\QQ$-Gorenstein situation.

\begin{proof}[Proof of Theorem \ref{thm.poschar}]
    We recall the description of the injective hull $E$ of the residue field in
    the case of a toric ring. It is the (graded) Matlis dual of $R$.
    Since $R$ consists of the monomials $x^m$ such that $(m,v_i)\geq 0$ for all $i$,
    its dual $E$ consist of the monomials $x^m$ such that $(m,v_i) \leq 0$.
    By \cite{WatFregFpure}, Theorem 2.5, $F^{*e}(E)$ is dual to
    $\omega^{(p^e-1)}_R$ where the latter consists of all $m \in M$ such that
    $(m,v_i) \geq 1-p^e$ for all $i$. Therefore,
    \[
        F^{*e}(E)=\langle x^m | (m,v_i)\leq p^e-1 \text{ for all $i$} \rangle.
    \]
    An easy consequence of Lemma \ref{lem.testIdeal} is that $x^m \in \tau(\aa^c)
    \iff x^{-m} \not\in 0^{*\aa^c}_E$. Clearly, since $x^m\cdot x^{-m}=1$ is nonzero in $E$
    the implication from left to right is clear. Conversely, if $x^m
    \not\in \tau(\aa^c) = \Ann_R{0^{*\aa^c}_{E_{R/m}}}$ there is $x^{m'} \in
    0^{*\aa^c}_{E_{R/m}}$ such that $x^{m+m'} \neq 0$ in $E_{R/m}$. Therefore,
    $x^{-(m+m')} \neq 0$ in $R$ and consequently, $x^{-m} = x^{-(m+m')}x^{m'}$ is in
    $0^{*\aa^c}_{E_{R/m}}$.
    Now the following chain of equivalences shows the result.
\[
    \begin{split}
    x^m \in \tau(\aa^c) &\iff x^{-m} \not\in 0^{*\aa^c}_E \iff \exists q :
    1\cdot x^{-qm}\aa^{\ceil{cq}} \neq 0 \text{ in }F^{*e}(E) \\
    &\iff \exists q : (-qm+\ceil{cq}\Newt(\aa)) \cap \{ m' | \forall i\,(m',v_i)\leq q-1 \}\cap M \neq
    \emptyset \\
    &\iff \exists q \exists w : qm+qw \in \ceil{cq}\Newt(\aa)\cap M \text{ and }\forall i\, (qw,v_i)
    \leq (q-1) \\
    &\iff \exists q,w : m+w \in \smallfrac{\ceil{cq}}{q}\Newt(\aa)\cap \smallfrac{1}{q}M \text{ and }\forall i\,
    (w,v_i)\leq 1-\smallfrac{1}{q} \\
    &\iff \exists w:\forall i\, (w,v_i) \leq 1 \text{ and } m+w \in \Int(\Newt(\aa^c)) \\
    \end{split}
\]
The second equivalence is just the definition of $\aa^c$ tight closure for the
zero submodule using that a toric ring is strongly $F$-regular, and thus 1 can
be used as a test element \cite[Theorem 1.7]{HaraYosh}. The only other
implication that needs explanation is the reverse of the last equivalence. For
this let $w$ be as in the last line. Since $m+w$ is in the interior we can
perturb $w$ slightly such that for large enough $q$ we have $qw \in M$,
$(w,v_i)\leq 1-\frac{1}{q}$ and $m+w \in (c+\frac{1}{q})\Newt(\aa) \subseteq
\frac{\ceil{cq}}{q}\Newt(\aa)$ as required.

Finally, $R$ is $\QQ$-Gorenstein if and only if there is $w_0$ such that for
all $v_i$ one has $(w_0,v_i)=1$. Then, clearly, the last condition is
equivalent to
\[
    m+w_0 \in \Int(\Newt(\aa^c))
\]
as claimed.
\end{proof}
Note that Theorems \ref{thm.ideal} and \ref{thm.poschar} also give a direct
argument in the toric case for the fact that $\tau(\aa^c)=\Jj(\aa^c)$ if $X$ is
$\QQ$-Gorenstein. More generally one can make the following observation, which
was pointed out by N. Hara:
\begin{corollary}
    Let $X$ be an affine toric variety over a field of positive characteristic and $\aa$ a monomial ideal. Then
    \[
        \tau(X,\aa^c) = \sum \Jj((X,\Delta),\aa^c)
    \]
    where the sum ranges over all effective torus invariant $\QQ$-divisors $\Delta$
    such that $K_X+\Delta$ is $\QQ$-Cartier.
\end{corollary}
\begin{proof}
    If $w \in M_\QQ$ and $\Delta$ are related by the formula $\Delta = \sum (1- (w,v_i))
    D_i$ we have that $K_X+\Delta = -\div x^w$
    is $\QQ$-Cartier and that $\Delta$ is effective if and only if $(w,v_i) \leq 1$
    for all $i$. By Theorem \ref{thm.poschar}, $x^m \in \tau(X,\aa^c)$ if and only if there is
    $w \in M_\QQ$ such that $(w,v_i) \leq 1$ and $m+w \in \Int(\Newt(\aa^c))$.
    With $\Delta$ as above this is equivalent to
    $x^m \in \Jj((X,\Delta),\aa^c)$ by Theorem \ref{thm.ideal}.
\end{proof}
Furthermore, it should be straightforward to define $\tau((X,\Delta),\aa^c)$ in
positive characteristic analogous to $\Jj((X,\Delta),\aa^c)$ in characteristic
zero. Then a test ideal version of Theorem 1 should hold. This would lead to
the statement, in analogy with the Corollary, that
\begin{equation}\label{eqn.taupair}
    \tau(X,\aa^c) = \sum_{\substack{\Delta \text{ effective,}\\ K_X+\Delta \text{ is $\QQ$-Cartier.}}}
    \tau((X,\Delta),\aa^c)
\end{equation}
for any toric variety $X$. This naturally leads to the following question: What
is the class of varieties such that equation (\ref{eqn.taupair}) holds?

On the other hand, in an effort to enlarge the definition of the multiplier
ideal beyond the $\QQ$-Gorenstein case, formula (\ref{eqn.taupair}) for the
test ideals (which does not depend on a $\QQ$--Gorenstein assumption) can serve
as a guideline. One is lead to speculate whether the following definition leads
to a useful generalization of the multiplier ideal, say if $X$ is over a field
of characteristic zero we set
\[
    \Jj(X,\aa^c) \defeq \sum_{\substack{\Delta \text{ effective,}\\ K_X+\Delta \text{ is $\QQ$-Cartier.}}} \Jj((X,\Delta),\aa^c)
\]
Clearly, in the $\QQ$-Gorenstein case this recovers the definition of
$\Jj((X,\emptyset),\aa^c)$ above. This definition would be justified, for
example, if $\Jj(X,\aa^c)$ satisfies a Nadel type vanishing theorem. At the
same time it is clear that one still needs to control the singularities of $X$,
at least one has to require that there exists an effective $\QQ$-divisor
$\Delta$ such that $K_X+\Delta$ is $\QQ$-Cartier in order for the definition to
make sense.

\bibliographystyle{amsalpha}

\bibliography{MultiplierToric}
\end{document}